\documentclass[letterpaper, 10 pt, conference]{ieeeconf}  

\IEEEoverridecommandlockouts                              

\usepackage{graphics} 
\usepackage{epsfig} 
\usepackage{mathptmx} 
\usepackage{times} 
\usepackage{amsmath} 
\usepackage{amssymb}  
\usepackage[ruled,vlined,linesnumbered]{algorithm2e}
\usepackage{url}
\usepackage{xcolor}
\usepackage{tikz}
\usetikzlibrary{arrows.meta, positioning}

\newtheorem{remark}{Remark}

\title{\LARGE \bf
Profit Maximization for a Robotics-as-a-Service Model}

\author{Joo Seung Lee$^{1}$ and Anil Aswani$^{1}$
\thanks{*This material is based upon work supported by the National Science Foundation under Grant No. DGE-2125913 and Grant No. CMMI-1847666.}
\thanks{$^{1}$Joo Seung Lee and Anil Aswani are with the department of Industrial Engineering and Operations Research,
        University of California, Berkeley, CA 94720
        {\tt\small \{jooseung\_lee,aaswani\}\@berkeley.edu}}%
}

\begin{document}

\maketitle
\thispagestyle{empty}
\pagestyle{empty}

\begin{abstract}

The growth of Robotics-as-a-Service (RaaS) presents new operational challenges, particularly in optimizing business decisions like pricing and equipment management. While much research focuses on the technical aspects of RaaS, the strategic business problems of joint pricing and replacement have been less explored. This paper addresses the problem of profit maximization for an RaaS operator operating a single robot at a time. We formulate a model where jobs arrive sequentially, and for each, the provider must decide on a price, which the customer can accept or reject. Upon job completion, the robot undergoes stochastic degradation, increasing its probability of failure in future tasks. The operator must then decide whether to replace the robot, balancing replacement costs against future revenue potential and holding costs. To solve this complex sequential decision-making problem, we develop a framework that integrates data-driven estimation techniques inspired by survival analysis and inverse optimization to learn models of customer behavior and robot failure. These models are used within a Markov decision process (MDP) framework to compute an optimal policy for joint pricing and replacement. Numerical experiments demonstrate the efficacy of our approach in maximizing profit by adaptively managing pricing and robot lifecycle decisions.
\end{abstract}

\section{Introduction}

Rapid advances in robotics and artificial intelligence (AI) are enabling the growth of \emph{cloud robotics} and \emph{robotics-as-a-service} (RaaS) \cite{chen2010robot,hu2012cloud,kehoe2015survey,buerkle2023towards}. The applications of RaaS to urban mobility \cite{straubinger2020overview,cohen2021urban} and drone delivery \cite{dorling2016vehicle,scott2017drone} are in particular finding rapid development and planned adoption, although concepts in agriculture \cite{dharmasena2019autonomous,milella2024robot}, manufacturing \cite{yan2017cloud}, and other domains \cite{hu2012cloud,toscano2022integrating} are also being considered. Yet despite growing attention to the RaaS paradigm, questions related to business aspects (e.g.,  pricing, maintenance, or replacement decisions) have been less well explored.

\subsection{Predictive Maintenance}

One important business aspect of RaaS is decisions related to maintenance. Because equipment maintenance is a major source of costs and delays in industry, there has been substantial interest in developing data-driven approaches for streamlining maintenance decisions. Predictive maintenance \cite{zhang2019data,carvalho2019systematic,wang2024dynamic} uses sensors that collect streaming data from equipment in order to predict degradation and identify maintenance that may be needed. Approaches based on more classical statistical modeling \cite{zhang2019data}, machine learning \cite{carvalho2019systematic}, and neural-network-based AI \cite{wang2024dynamic} have been proposed. Such approaches have found high-visibility success in the aviation context, where predictive maintenance is used for commercial aircraft and their constituent components \cite{ge_maintenance_insight}.

\subsection{Failure Modeling}

Another important business aspect of RaaS is decisions related to replacement of robots. Whereas questions of maintenance are related to understanding tradeoffs between the timing and amount of maintenance with uncertain knowledge about the level of equipment degradation, decisions regarding replacement are related to managing tradeoffs between holding costs (e.g., costs related to storage), failure costs (e.g., costs related to failures while in operation), and replacement costs (i.e., costs for purchasing new robots). 

Though there is relatively little work on managing replacements, there has been considerably more research on modeling rates of survival or failure. Component-based models for predicting failure rates in electronics have been developed \cite{mil-hdbk-217f,hinton2006prism} and found use in practice. Survival analysis \cite{vittinghoff2012regression,karim2019reliability} uses data to empirically estimate rates of failure and its sensitivity to predictive factors, and it is most commonly used within healthcare \cite{vittinghoff2012regression,collett2015modelling,eloranta2021cancer} with less frequent use in other domains like predicting machine failure \cite{karim2019reliability,papathanasiou2023machine}. The mainstay approach for survival analysis is the proportional hazards model \cite{cox1972regression,kumar1994proportional}, which is able to successfully handle the right-censoring of survival data (i.e.,  life-span data is not available for items or people that have survived at the time of data collection).

\subsection{Inverse Optimization}

An important business aspect of RaaS is decisions related to pricing. Imagining a fleet of robots that perform different jobs arriving sequentially, a relevant business decision is on how to price each job in order to maximize profit while considering equipment degradation and other costs. The most related line of work is the topic of data-driven contract design. In the area of contract design \cite{laffont2002theory}, the \emph{principal} is a leader in making decisions while the \emph{agent} is a follower that makes decisions knowing the decision of the principal. The principal seeks to make an optimal decision knowing the \emph{leader-follower} information structure.

One set of approaches to data-driven contract design make use of inverse optimization \cite{aswani2018inverse,chan2025inverse}, which are techniques that input decision data from agents and estimate utility function models for which the decision data are optimal. These estimated utility function models can then be input into bilevel optimization problems (i.e., optimization problems in which some constraints are optimization problems themselves) in order to design an optimal contract knowing the leader-follower information structure \cite{aswani2019data,mintz2023behavioral,li2025adaptive}.

Another approach to data-driven contract design revolves around multi-armed bandit models \cite{dogan2023estimating,dogan2023repeated,scheid2024incentivized} in order to better understand exploration-exploitation issues. Based upon multi-armed bandit models in which choosing an action provides a stochastic reward, in these models the agent chooses an action in order to maximize their own expected rewards while considering an incentive amount provided by the principal for each action. The central issue is in how the principal can indirectly estimate the agent's expected rewards for each action having observations of only the sequence of actions chosen and not directly measuring the agent's realized rewards. Viewed from a certain perspective, this estimation problem is analogous to inverse optimization. 

\subsection{Contribution and Outline}

To the best of our knowledge, this is the first paper to consider the problem of maximizing profit for robots operating as part of a RaaS business where pricing and replacement decisions are made. Potential customers arrive sequentially, where each customer's job has a vector-valued context and a model for (a) the maximum price willing to pay for the job and (b) the stochastic degradation of a robot completing the job. The operator of the RaaS business must decide what price to charge for servicing the job, and the customer gets to decide whether or not to accept the price. If a customer rejects a price, then the job is lost forever. If a customer accepts a price, then a robot services the job. The robot fails during the job with a stochastic rate depending upon its current degradation level. If the robot completes the job, then the operator gets to decide whether or not to replace the robot considering the tradeoffs between robot replacement costs, holding (e.g., storage) costs, and opportunity costs (i.e., possible future revenue without replacing the robot). 

The paper is outlined as follows: Section \ref{sec:s&m} formally defines the elements of the quantitative model that we consider in this paper, and it presents the profit maximization problem formulated as a Markov decision process (MDP). Section \ref{sec:sf} describes our proposed approaches for optimizing the RaaS operator's decisions within the MDP. The paper concludes with numerical experiments in Section \ref{sec:ne} that demonstrate the potential efficacy of our proposed approach.

\section{Model and Problem Statement}
\label{sec:s&m}

We consider an operator managing an RaaS business in an online setting that operates with a single robot at a time. This section presents a model for customers and their jobs as well as a model for the operator of the RaaS. The section concludes by formally defining the goal of developing a policy that the RaaS operator can use for pricing and robot replacement that maximizes profit in the face of uncertainty about customer preferences and robot degradation.

\subsection{Customer and Jobs Model}

In our online setting, customers arrive sequentially at times indexed by $k=1, 2, \dots$. The $k$-th customer is characterized by a tuple $(x_k, T_k, \tau_k)$, where:
\begin{itemize}
    \item $x_k \in \mathbb{R}_+^d$ is a context vector representing characteristics of the customer's job, with a normalization of $\|x_k\|_2 \leq 1$;
    \item $T_k > 0$ is the desired rental duration for the job;
    \item $\tau_k > 0$ is the interarrival time since the last customer interaction.
\end{itemize}
Upon arrival of a customer, the RaaS operator observes $(x_k, T_k)$ and offers a price $p_k$ to the customer.

We assume that a customer's decision to rent a robot for their job is deterministic and governed by an unknown (to the RaaS operator) utility vector $u \in \mathbb{R}_+^d$, with $\|u\|_2 \leq 1$. The customer accepts the offered price if their utility exceeds the price $u^\top x_k \geq p_k$. One of the RaaS operator's challenges is learning $u$ in order to be able to set profit-maximizing prices.

\begin{remark}
    There are two key aspects to the above assumption. The first is that we are assuming that a customer's utility for their job is a linear function of the context $x_k$. Such a linear context assumption is common in the literature involving online models \cite{agrawal2013thompson,dimakopoulou2019balanced,bastani2021mostly}. 
\end{remark}

\begin{figure*}
    \centering
    \includegraphics[width=\linewidth]{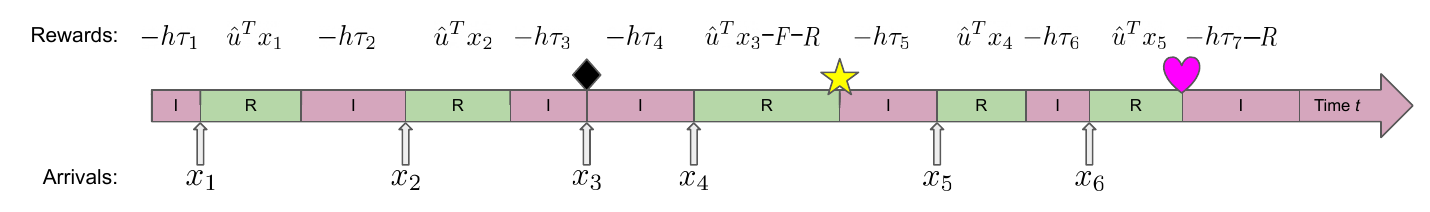}
    \caption{An example timeline illustrating the operational flow for the Robotics-as-a-Service (RaaS) provider. Upward arrows indicate the arrival of new customers ($x_i$). The robot's status alternates between being under rental ('R' in a green box) and idle ('I' in a purple box). The timeline also depicts key events like a price rejection from an arriving customer (diamond symbol), an unexpected robot failure during a job (star symbol), and a voluntary replacement initiated by the operator (heart symbol). The rewards associated with each period are shown above the timeline, reflecting revenues $\hat u^Tx_i$, holding costs $h\tau_i$, and costs for failure $F$ and replacement $R$.}
    \label{fig:timeflow}
\end{figure*}

\begin{remark}
The second is that we are assuming a customer rents a robot as long as their utility for completing the job is greater than \emph{or equal to} the price. In the case where the utility equals the price, presumably the customer receives no surplus value for the completion of the job since they must pay an amount equal to their utility. However, we note that such an assumption is common in the contract design literature \cite{laffont2002theory}, and the corresponding inequality is often referred to as a \emph{participation constraint}. We also note that this assumption can be relaxed to an alternative model where a customer rents a robot if $u^\top x_k \geq p_k + \nu$ where $\nu$ is a constant that describes the minimum surplus value a customer must receive for accepting a rental, but for conciseness we assume the above.
\end{remark}

\subsection{Robot Degradation and Failure Model}
A robot's condition is modeled by its cumulative operational history. Let $X$ be the sum of context vectors from all prior successful rentals since the last replacement, and let $t_{\text{age}}$ be the sum of rental durations from all prior successful rentals since the last replacement. During a new rental with context $x_k$, the robot is subject to breakdown. We model this using a Cox proportional hazards model \cite{cox1972regression,kumar1994proportional}: The \emph{hazard rate} at time $t \in [0, T_k]$ into the rental, which can be interpreted as the probability the robot will not survive in an infinitely small period of time beyond $t_{\text{age}} + t$ conditioned on the robot having survived up to time $t_{\text{age}} + t$, is
\begin{equation}
    \lambda(t | X, x_k) = \lambda_0(t + t_{\text{age}}) \exp(\theta^\top (X + x_k))
\end{equation}
where
\begin{itemize}
    \item $\lambda_0(\cdot)$ is an unknown (to the RaaS operator), non-negative baseline hazard function;
    \item $\theta \in \mathbb{R}_+^d$ is an unknown (to the RaaS operator) degradation parameter vector, with normalization $\|\theta\|_\infty \leq 1$, that modulates the impact of cumulative context on failures.
\end{itemize}
A breakdown incurs a failure cost $F$ and a replacement cost $R$, after which the robot's cumulative context and cumulative rental duration are reset to states corresponding to a new robot (i.e., $X=0$ and $t_{\text{age}} = 0$). The second challenge for the RaaS operator is to learn $\theta$ and $\lambda_0(\cdot)$ to effectively manage the risk of breakdown.

\begin{remark}
    Note that we are assuming that the hazard rate depends upon a linear function of the cumulative contexts $X$, which is the sum of the contexts from prior successful rentals since the last replacement. This is a modest generalization of the linear context assumption that is common in the literature involving online models \cite{agrawal2013thompson,dimakopoulou2019balanced,bastani2021mostly}, but it is arguably reasonable since we would expect failure rates to depend upon the history of completed jobs.
    \end{remark}

\subsection{Model of RaaS Operator}
\label{sec:mraas}

We formulate the RaaS operator's decision problem as an infinite-horizon  MDP that is defined by the tuple $(\mathcal{S}, \mathcal{A}, \mathcal{P}, \mathcal{R}, \gamma)$.

\textbf{State Space $\mathcal{S}$} The state $s \in \mathcal{S}$ must capture the robot's condition and the current decision epoch. We define $s = (X, t_{\text{age}}, \phi)$, where:
\begin{itemize}
    \item $X \in \mathbb{R}_+^d$ is the cumulative context vector since the last replacement;
    \item $t_{\text{age}} \in \mathbb{R}_+$ is the cumulative operational time since the last replacement, conditioning the baseline hazard;
    \item $\phi \in \{0, 1\}$ is a phase indicator; $\phi=0$ denotes a customer has just arrived (i.e., pricing decision), and $\phi=1$ denotes the robot is idle immediately after a rental termination and the operator has to make a choice regarding a replacment decision.
\end{itemize}

\textbf{Action Space $\mathcal{A}(s)$} The actions available to the RaaS operator depend on the phase $\phi$, which is categorized into:
\begin{itemize}
    \item Case $\phi=0$ (i.e., arrival): Given an arriving customer with context $x_k$, the action is $a \in \{\text{Accept}, \text{Reject}\}$;
    \item Case $\phi=1$ (i.e., idle): The action is $a \in \{\text{Continue}, \text{Replace}\}$. We assume `Replace' incurs cost $R$ and resets the state.
\end{itemize}

\textbf{Reward Function $\mathcal{R}$ \& Transition Dynamics $\mathcal{P}$} Let the current state be $s=(X, t_{\text{age}}, \phi)$. The rewards and state transitions are determined by the system's current state $s=(X, t_{\text{age}}, \phi)$ and the principal's chosen action. As we will see, to maximize long-run profit per unit time, we normalize each immediate reward by dividing it by the elapsed time in the transition (rental duration or interarrival time).

When the system is in the arrival phase (i.e., Case $\phi=0$), a customer with characteristics $(x_k, T_k)$ has just arrived. If the RaaS operator chooses to `Reject' the rental, the system transitions to $(X, t_{\text{age}}, 0)$ with zero reward and elapsed time.

\begin{figure*}[!htb]
    \centering
    \includegraphics[width=0.9\linewidth]{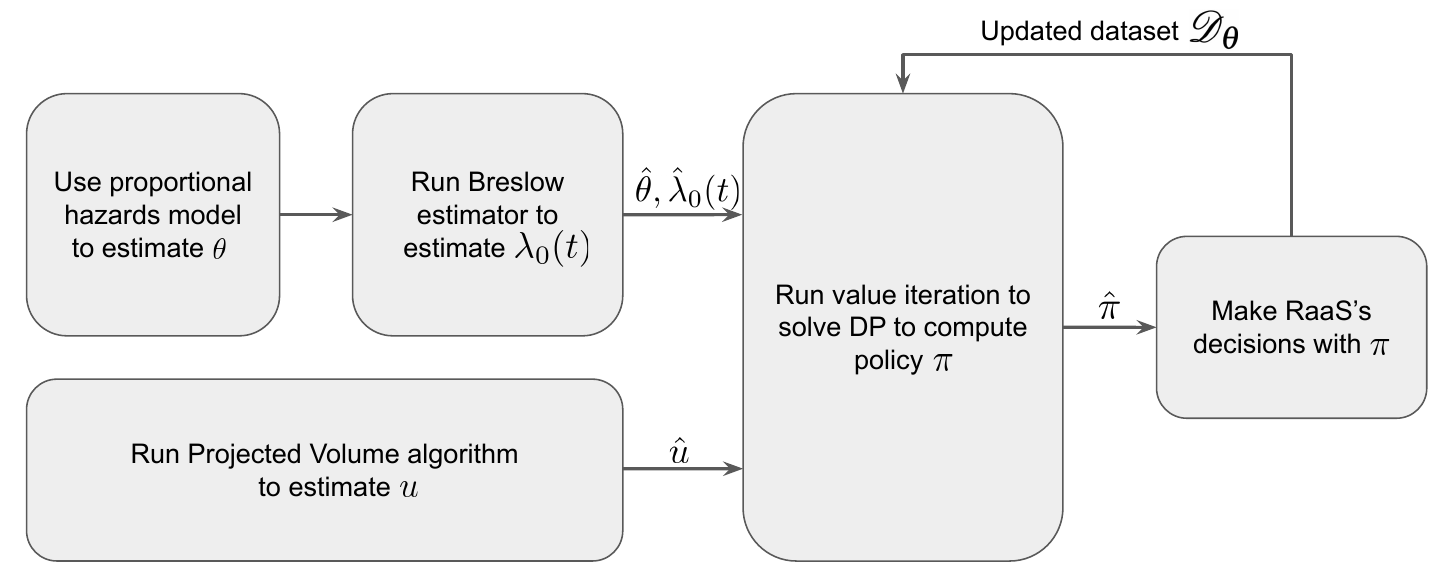}
    \caption{A diagram of the three-phase learning and control framework. \textbf{Phase 1:} The Projected Volume algorithm is run to obtain a reliable estimate of the customer utility vector, $\hat u$. \textbf{Phase 2:} This estimate is used to learn the degradation parameters $(\hat\theta, \hat\lambda_0)$ and train an initial control policy, $\hat\pi$. \textbf{Phase 3:} The system operates using this policy until a predefined number of new robot failures ($K$) occur. The data from these failures updates the dataset $\mathcal{D}_\theta$, which is then used to re-estimate the parameters and retrain the policy in a repeating cycle.}
    \label{fig:phases}
\end{figure*}

If the RaaS operator chooses to `Accept', the customer makes an acceptance decision. If $p_k > u^\top x_k$ (i.e., the price is greater than the customer's utility for the job), then the customer rejects the offered price and the system transitions to $(X, t_{\text{age}}, 0)$ with zero reward and elapsed time. Otherwise, if  $p_k \leq u^\top x_k$ then the customer accepts the offered price, and the rental begins. The outcome is stochastic, and it is specifically governed by the probability of robot failure $p_{\text{fail}} = 1 - \exp(-\exp({\theta}^\top(X+x_k)) \int_0^{T_k} {\lambda}_0(t_{\text{age}}+t) dt)$. With probability $1-p_{\text{fail}}$, the rental is successful, yielding a reward of $p_k/T_k$ (i.e., the rental revenue). The system state is updated to reflect the accumulated usage, transitioning to the idle phase at state $(X+x_k, t_{\text{age}}+T_k, 1)$. With probability $p_{\text{fail}}$, the robot breaks down during the rental, with time to failure $f_k < T_k$. This results in a reward of $(p_k - F - R)/f_k$, accounting for revenue, failure cost, and replacement cost. The robot is replaced, and the state resets to $(0, 0, 1)$ because the RaaS operator is now using a new robot.

When the system is in the idle phase (i.e., Case $\phi=1$), the RaaS operator decides whether or not to replace the robot. If the action is to `Continue', then no replacement occurs. The system incurs a holding cost over the interarrival time $\tau_{k+1}$, for a reward of $-h\tau_{k+1}/\tau_{k+1}=-h$, and transitions to state $(X, t_{\text{age}}, 0)$ when the next customer arrives. If the action is to `Replace', then the RaaS operator pays the replacement cost $R$ in addition to the holding cost. The reward is $-R/\tau_{k+1} - h$, and the state is reset to $(0, 0, 0)$ for the next customer arrival. 


\subsection{Problem Statement}
We assume that the RaaS operator's objective is to maximize their average profit per unit of time. Formally, let $\pi$ denote a policy that specifies, for each arriving customer $k$, whether to offer a rental (and at what price $p_k$) and, after each rental or rejection, whether to voluntarily replace the robot. Under policy $\pi$, the system generates a sequence of revenues and costs over an infinite horizon. Fig. \ref{fig:timeflow} provides an example timeline of this process, illustrating the sequence of events and their associated rewards.  Let $\mathcal{T}[N]$ denote the total time elapsed after the first $N$ customers:
\begin{equation}
\mathcal{T}[N] = \sum_{k=1}^N \tau_k + \sum_{k=1}^N T_k \cdot \mathbb{I}\{a_k = \text{accept}\},
\end{equation}
accounting for both idle and rental periods (adjusted for partial rentals due to failures). The revenue consists of accepted prices $\{p_k : a_k = \text{accept}\}$, where $a_k$ is the customer's acceptance decision. The costs include failure penalties $F$ for each mid-rental breakdown, replacement costs $R$ for each breakdown or voluntary replacement, and holding costs $h \cdot \tau_k$ for each interarrival idle period $\tau_k$.

The fundamental optimization problem that the RaaS operator seeks to solve is thus given by
\begin{multline}
\label{eqn:raasobj}
    \max_{\pi} \lim_{N \to \infty} \frac{1}{\mathcal{T}[N]} \Bigg[ \sum_{k=1}^N p_k \cdot \mathbb{I}\{a_k = \text{accept}\} +\\
    - F \cdot n_F[N] - R \cdot n_R[N] - h \sum_{k=1}^N \tau_k \Bigg],
\end{multline}
where $n_F[N]$ is the number of failures and $n_R[N]$ is the number of replacements (including those triggered by failures) up to $N$, and $\mathbb{I}\{\cdot\}$ is the indicator function. All terms depend on $\pi$ by its influence on acceptances, failures, and replacements. Since the unknown parameters $(u, \theta, \lambda_0)$ must first be learned, the principal has to manage the exploration-exploitation tradeoff between balancing exploration to estimate these parameters with exploitation to maximize profit. 

\section{Solution Framework}
\label{sec:sf}

To tractably compute a performant policy that is able to optimize the RaaS operator's objective \ref{eqn:raasobj}, we propose an algorithm that decouples the optimization problem into three sequential phases: (1) estimation of the customer utility vector, (2) estimation of the robot degradation parameters, and (3) computation of an optimal control policy. This iterative learning and control process is illustrated in Fig. \ref{fig:phases}. We manage the exploration-exploitation tradeoff by making use of a standard $\epsilon$-greedy approach \cite{dann2022guarantees}. This section describes the algorithms used for each phase, and implementation details are provided in the next section.

\subsection{Phase 1: Learning the Utility Parameters $u$}
To learn the utility parameters $u$, we employ an active learning strategy based on the projected volume algorithm \cite{scheid2024incentivized, lobel2018multidimensional}. This method treats the estimation as a multidimensional binary search problem, maintaining a convex uncertainty set $S_k \subset \{s \in \mathbb{R}^d : |s|_2 \leq 1\}$ guaranteed to contain the true vector $u$. It also tracks a set of orthogonal directions where the uncertainty set has small width. At each customer arrival $k$ with context $x_k$, the algorithm first ``cylindrifies" the set by expanding it along these small-width directions (effectively projecting onto the subspace of large uncertainty to form a high-dimensional cylinder). It then computes an approximate centroid $\widehat{s}_k$ of the cylindrified set via sampling and offers a query price $p_k = \hat s_k^\top x_k$. The customer's accept/reject decision provides a linear constraint ($u^\top x_k \geq p_k$ or $u^\top x_k < p_k$), which prunes the uncertainty set via a volume-cutting halfspace: $S_{k+1} = S_k \cap \{s : s^\top x_k \gtrless p_k\}$. This phase continues until the uncertainty set is sufficiently small, yielding the final centroid as a high-confidence estimate $\widehat{u}$.

\subsection{Phase 2: Learning Degradation Parameters $(\theta, \lambda_0(\cdot))$}
Once a reliable estimate $\widehat{u}$ is obtained, we enter a data collection phase to learn the degradation model. Using a fixed pricing policy (e.g., $p_k = \widehat{u}^\top x_k$), we gather rental data tuples of (cumulative context, rental duration, cumulative operational time, breakdown indicator).

\textbf{Estimating ${\theta}$}: We estimate the degradation parameter $\theta$ by maximizing the partial log-likelihood of the Cox model \cite{vittinghoff2012regression,cox1972regression}:
\begin{multline}
    \mathcal{L}(\theta) = \sum_{t\in \mathcal{F}}\Bigg[ \sum_{k\in\mathcal{D}(t)} \theta^\top(X_k + x_k) + \\
    -|\mathcal{D}(t)|\log \Big( \sum_{j \in R(t)} \exp(\theta^\top(X_j + x_j))\Big)\Bigg],
\end{multline}
where $\mathcal{F}$ is the set of unique failure times, $\mathcal{D}(t)$ is the set of failures at time $t$, $R(t)$ is the risk set $t$ (rentals active at $t$), and $\mathcal{D}(t)$ is the number of failures at $t$. This conveniently eliminates the dependency on the baseline hazard $\lambda_0(t)$. We impose the constraint $\theta \geq 0$ to ensure that usage is strictly degrading under estimated $\widehat\theta$. 

\textbf{Estimating ${\lambda}_0(\cdot)$}: With the estimate $\widehat{\theta}$, we then estimate the cumulative baseline hazard $\Lambda_0(t) = \int_0^t \lambda_0(s) ds$ using the non-parametric Breslow estimator \cite{breslow1972contribution}. The instantaneous baseline hazard $\widehat{\lambda}_0(t)$ is recovered by smoothing the increments of $\widehat{\Lambda}_0(t)$ and then differentiating.

\subsection{Phase 3: Solving the MDP Problem}

If the parameters $(u,\theta,\lambda_0(\cdot))$ are known, then in principle the optimal policy can be computed using dynamic programming \cite{Bertsekas2017VolI,Bertsekas2012VolII}. The relevant computational issue is whether or not the dimensionality of the state space is small enough to employ an algorithm that directly computes an optimal policy via dynamic programming.

Though the state space (see Section \ref{sec:mraas})  is continuous and high-dimensional (because of the cumulative context vector $X\in\mathbb{R}^d_+$), there exists a lower-dimensional projection that preserves all information relevant for decision-making. The key insight is that the rewards and dynamics depend only on scalar projections of the context vectors. We thus redefine the state for the arrival phase (i.e., Case $\phi=0$) as a 4-dimensional tuple $s' = (c_c + c_x, c_u, T, t)$, where:
\begin{itemize}
    \item $c_c = \theta^\top X$ represents the cumulative degradation context;
    \item $c_x = \theta^\top x_k$ is the customer-specific degradation context;
    \item $c_u = u^\top x_k$ is the customer-specific revenue context;
    \item $T$ is the job's rental duration;
    \item $t$ is the robot's cumulative active lifetime.
\end{itemize}
For the idle phase (i.e., Case $\phi=1$), the state is fully described by the cumulative degradation $c_c$ and $t$. This projected state space is small enough that dynamic programming can be applied computationally in a direct manner.

\begin{figure*}
    \centering
    \includegraphics[width=\linewidth]{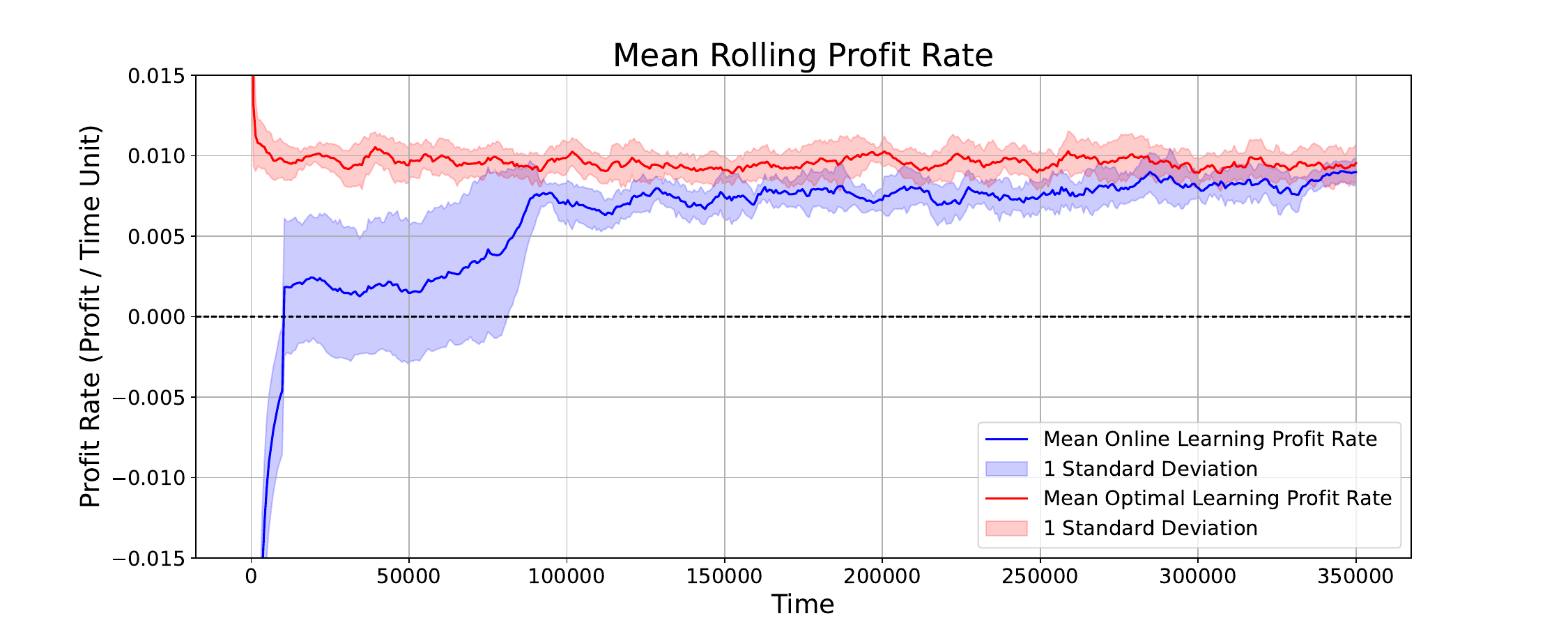}
    \caption{Rolling average profit rate per time unit for the online learning policy and the oracle optimal policy, computed over a 10,000-unit window and averaged across 10 simulation runs with shaded standard deviations.}
    \label{fig:profit_rate_comp}
\end{figure*}

Our approach to computing a policy for the MDP, which optimizes the RaaS operator's objective (\ref{eqn:raasobj}), is informed by the above state space projection. Specifically, we use the approaches described in Phases 1 and 2 to estimate the parameters $(u,\theta,\lambda_0(\cdot))$. Value iteration  \cite{Bertsekas2012VolII}, over a finite, discretized grid, is used to compute an optimal policy for the MDP with model parameters set to be equal to the estimated parameters $(\widehat u, \widehat\theta, \widehat\lambda_0)$. The flow between the various phases of our solution framework are summarized in Fig. \ref{fig:phases}.

\begin{algorithm}
\SetAlgoLined
\DontPrintSemicolon
\SetKwInput{KwIn}{Input}
\SetKwInput{KwOut}{Output}
\SetKwFunction{FUpdateU}{ProjectedVolumeUpdate}
\SetKwFunction{FEstTheta}{EstimateDegradationParameters}
\SetKwFunction{FTrainPolicy}{TrainPolicyFQI}
\SetKwFunction{FExecPolicy}{ExecutePolicy}

\KwIn{
    Total customers $N$, Termination thresholds $\epsilon_u$, Policy update threshold $K$, Customer distribution $\mathcal{D}_{cust}$, Costs $h$, $F$, $R$, Initial datasets $\mathcal{D}_u \leftarrow \emptyset$, $\mathcal{D}_\theta \leftarrow \emptyset$\;
}
\KwOut{History $\mathcal{H}$}

Initialize $S_0 \leftarrow \{s \in \mathbb{R}^d : \|s\|_2 \leq 1\}$, $V_0 \leftarrow \emptyset$\;
$\widehat{u} \leftarrow \mathbf{0}$, $\widehat{\theta} \leftarrow \mathbf{0}$\;
$M_{state} \leftarrow (X=\mathbf{0}, t_{age}=0)$\;
$\mathcal{H} \leftarrow \emptyset$\;
utility\_learning\_done $\leftarrow$ False\;

\For{$k \leftarrow 1$ \KwTo $N$}{
    Observe customer $(x_k, T_k, \tau_k) \sim \mathcal{D}_{cust}$\;
    \If{not utility\_learning\_done}{
        Take centroid $\widehat{u}$ of $S_{k-1}$\;
        Offer $p_k = \widehat{u}^\top x_k$, observe $a_k \in \{\text{accept}, \text{reject}\}$\;
        Append $(x_k, p_k, a_k)$ to $\mathcal{D}_u$\;
        Do Projected Volume Update to get $S_k, V_k$\;
        \If{$a_k = \text{accept}$}{
            Observe outcome $o_k \in \{\text{success}, \text{failure}\}$\;
            Append $(M_{state}, T_k, o_k)$ to $\mathcal{D}_\theta$\;
            Update $M_{state}$\;
        }
        \If{$\text{diam}(S_k, x_k / \|x_k\|_2) < \epsilon_u$}{
            Estimate $\widehat{\theta}, \widehat{\lambda}_0$ from $\mathcal{D}_\theta$\;
            Train policy $\pi$ using $\widehat{\theta}, \widehat{\lambda}_0, \widehat u$\;
            Mark utility learning complete\;
        }
    }\Else{
        Query action $a_{k,arr}$ from $M_{state}, x_k, T_k$ from $\pi$\;
        \If{$a_{k,arr} = \text{Accept}$}{
            Offer $p_k = \widehat{u}^\top x_k$, observe $a_k$\;
            \If{$a_k = \text{accept}$}{
                Observe $o_k$\;
                Append $(M_{state}, T_k, o_k)$ to $\mathcal{D}_\theta$\;
                Update $M_{state}$\;
            }
        }
        Query action $a_{k,dep}$ from $M_{state}$ from $\pi$\;
        \If{$a_{k,dep} = \text{Replace}$}{
            Reset $M_{state}$\;
        }
        \If{$K$ failures since last $\pi$ update}{
            Estimate $\widehat{\theta}, \widehat{\lambda}_0$ from $\mathcal{D}_\theta$\;
            Train policy $\pi$ using $\widehat{\theta}, \widehat{\lambda}_0, \widehat u$\;
        }
    }
    Append to $\mathcal{H}$\;
}
\Return{$\mathcal{H}$}\;
\caption{RaaS Operator's Algorithm}
\label{alg}
\end{algorithm}

\section{Numerical Experiments}
\label{sec:ne}
\subsection{Experimental Setup}
Our empirical validation is performed within a discrete-time simulation environment designed to validate our online learning and control framework coded using the Python 3 programming language and run on an Apple M4 Pro.

\begin{figure}
    \centering
    \includegraphics[width=\linewidth]{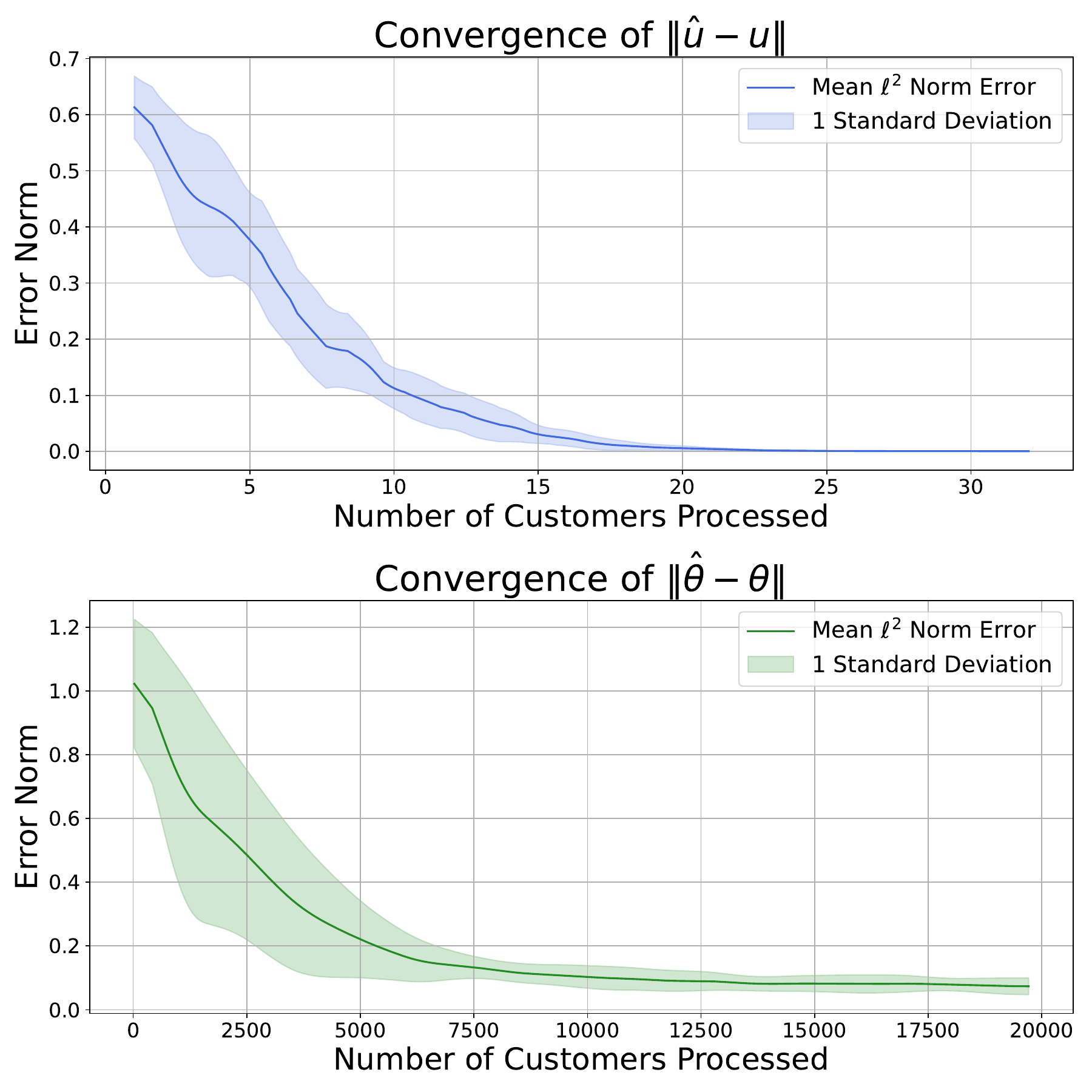}
    \caption{$\ell_2$-norm estimation errors for utility vector $\widehat{u}$ (top) and degradation parameter $\widehat{\theta}$ (bottom) versus number of customers processed, showing means and standard deviations over 10 runs.}
    \label{fig:convergence}
\end{figure}

We simulate a horizon of $N=20000$ sequential customer arrivals. The customer context vectors $x \in \mathbb{R}^d$ are drawn from a uniform distribution over the positive orthant of a $d$-dimensional unit ball, with the dimensionality set to $d=4$. Customer interarrival times $\tau_k$ and their desired rental durations $T_k$ are stochastic, sampled from exponential distributions with mean values of $5$ and $10$ time units, respectively. The financial parameters of the simulation are set as follows: the holding cost per unit time is $h=0.02$, the cost incurred upon a robot failure is $F=0.75$, and the cost for a voluntary replacement is $R=1.5$. The true customer utility vector is $u=(0.37, 0.11, 0.34, 0.71)$. The robot's failure mechanism is governed by a Cox proportional hazards model with a true degradation parameter given by $\theta =(0.5,0.2,0.4,0.3)$ and by a true baseline hazard $\lambda_0(t)=\lambda$ modeled by exponential process with a constant rate $\lambda=0.001$.

\subsection{Learning Scheme}
We compute a policy in our simulation using Algorithm \ref{alg}, which combines the three phases of the solution framework from Section \ref{sec:sf}. Experiment-specific details (to support reproducibility and address practical challenges) are as follows:

\textbf{Phase 1: Utility Vector Estimation.} Building on the Projected Volume algorithm, we allocate an initial exploration period that typically spans the first $\sim$40 customer arrivals, terminating dynamically when the diameter of the uncertainty set $S_k$ (measured along the normalized context direction $x_k / \|x_k\|_2$) falls below $10^{-4}$. The centroid of $S_k$ is approximated using hit-and-run sampling with 2000 samples, a burn-in of 500$\cdot d^2$ steps, and a thinning factor adjusted for efficiency (tolerance $10^{-4}$, target acceptance rate 0.01). During this phase, any accepted rentals contribute initial degradation data to $\mathcal{D}_\theta$ if failures occur.

\begin{figure*}
    \centering
    \includegraphics[width=\linewidth]{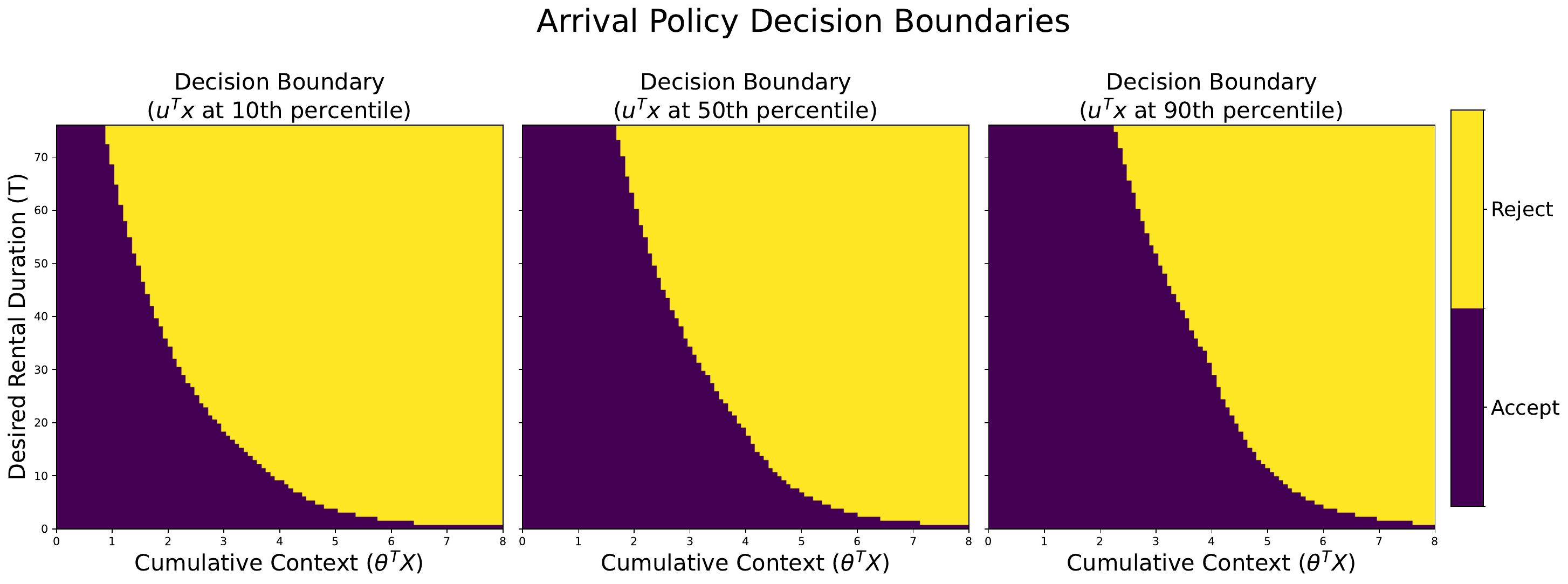}
    \caption{Decision regions for job acceptance (purple) or shutdown (yellow) in the arrival phase, plotted against desired rental duration $T$ and cumulative degradation $\theta^\top X + \theta^\top x$, at selected percentiles of customer revenue $u^\top x$.}
    \label{fig:arrival_policy}
\end{figure*}

\textbf{Interleaved Phases 2 and 3: Degradation Estimation and Policy Learning.} Immediately after Phase 1 convergence, we estimate $\widehat{\theta}$ by maximizing the partial log-likelihood under the constraint $\theta \geq 0$. For the baseline hazard, we apply the Breslow estimator to compute $\widehat{\Lambda}_0(t)$, followed by smoothing its increments via a Gaussian kernel density estimator to obtain a continuous $\widehat{\lambda}_0(t)$. If insufficient failures have been observed (fewer than 2), we continue data collection using near-optimal pricing $p_k = \widehat{u}^\top x_k - \varepsilon_p$ ($\varepsilon_p=10^{-2}$) to encourage accepts until the threshold is met, then proceed to estimating $\theta$ and $\lambda_0(\cdot)$.

With the system dynamics and reward models specified by the learned parameters $(\widehat u, \widehat\theta, \widehat\lambda_0)$, we compute the optimal policy by solving the MDP directly. This is accomplished using value iteration on a discretized representation of the four-dimensional state space $(c_c+ c_x, c_u, T, t)$. To manage the uncertainty from initial parameter estimates, the agent begins with a decaying $\epsilon$-greedy policy derived from the initial value function with initial $\epsilon = 0.10$ and decay factor $\gamma_d=0.95$. This encourages exploration and mitigates the risk of a greedy policy exploiting inaccuracies in the early estimates  $\widehat\theta$ and $\widehat\lambda_0(\cdot)$, which might otherwise lead to suboptimal, premature `Reject' actions. The policy runs until $K=100$ additional failures accumulate, at which point we re-estimate $\widehat{\theta}$ and $\widehat{\lambda}_0(\cdot)$ from the updated $\mathcal{D}_\theta$ and retrain the policy.  This batch update process repeats throughout the horizon, adapting to new degradation data.

\subsection{Results}
\subsubsection{Profitability}
To benchmark our learning algorithms, we also train an oracle optimal policy which assumes \emph{perfect} knowledge of $u, \theta$, and $\lambda_0(t)$ and run simulation under the identical environment, representing an upper bound on performance. We evaluate the profitability of our online learning framework by comparing its rolling average profit rate (with window size 10,000) to that of the oracle optimal policy with full knowledge of the true parameters $u, \theta$, and $\lambda_0(t)$. Results are aggregated over 10 independent simulation runs, with means and standard deviations reported. As shown in Fig. \ref{fig:profit_rate_comp}, the optimal policy maintains a consistent positive profit rate from the outset, reflecting its perfect alignment with the environment. In contrast, the online policy initially yields negative rates due to exploratory pricing for estimating $u$ and suboptimal decisions from early, inaccurate estimates of $\theta$. However, as learning progresses through batch updates, the online policy's profit rate recovers, becoming positive and converging toward the optimal rate. This demonstrates the framework's ability to adaptively achieve near-optimal performance despite initial uncertainty.

\subsubsection{Convergence of $\widehat\theta$ and $\widehat u$}
To assess parameter estimation, we track the $\ell_2$-norm error between the estimates $\widehat{u}$ and $\widehat{\theta}$ and their ground-truth values, averaged over 10 runs with standard deviations. Fig. \ref{fig:convergence} shows rapid convergence for $\widehat{u}$, reaching near-zero error after approximately 40 customers, satisfying the termination criterion for Phase 1. Convergence for $\widehat{\theta}$ is slower and asymptotic, reflecting the need for accumulating failure observations. Though estimation error decreases slowly after an initial rapid decrease, Fig. \ref{fig:profit_rate_comp} indicates that this does not substantially harm the profit rate.

\begin{figure}
    \centering
    \includegraphics[width=\linewidth]{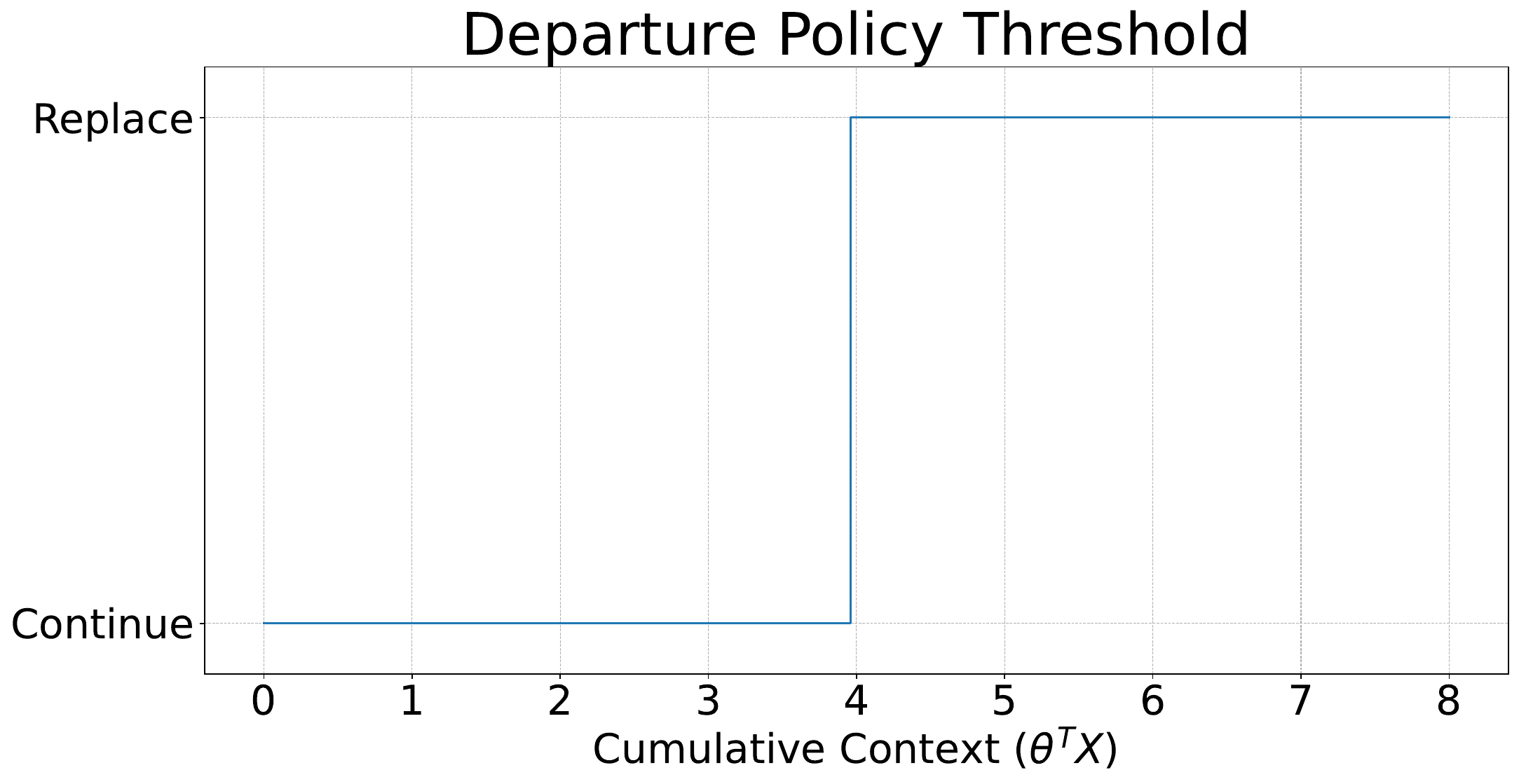}
    \caption{Replacement threshold in the departure phase as a function of cumulative degradation $\theta^\top X$, with `Replace' above and `Continue' below the boundary.}
    \label{fig:departure_policy}
\end{figure}

\subsubsection{Policy Structure}
The computed optimal policy exhibits intuitive structures in its arrival and departure phases, balancing revenue potential against degradation risks, failure costs, holding costs, and replacement expenses.

In the arrival phase (see Fig. \ref{fig:arrival_policy}), decision boundaries for accepting or shutting down jobs are visualized in the plane of desired rental duration $T$ and projected cumulative degradation $\theta^\top (X+x)$, with slices at the 10th, 50th, and 90th percentiles of projected customer revenue $u^\top x$. Higher revenue percentiles expand the acceptance region, as increased expected payoffs offset degradation-related risks for a given failure probability (determined by $T$ and $\theta^\top(X+x)$).

In the departure phase (see Fig. \ref{fig:departure_policy}), a threshold at $\theta^\top X \approx 4$ governs replacement decisions. Below this, retaining the robot is optimal to capitalize on potential future high-revenue jobs; above it, replacement is preferred to preempt likely failures (incurring additional costs beyond replacement) or excessive holding costs while awaiting suitable customers. 

\section{Conclusion}
Our framework enables RaaS operators to learn customer utilities and robot degradation models online, yielding policies that converge to near-optimal profitability through adaptive pricing and replacement decisions. Numerical results validate its efficacy in balancing revenue against operational risks. Future extensions could consider more realistic settings such as stochastic customer acceptance behaviors (see for example \cite{dogan2023estimating}), which necessitates more complex models.





\bibliography{bibliography}

\begin{thebibliography}{10}

\bibitem{chen2010robot}
Y.~Chen, Z.~Du, and M.~Garcia-Acosta, ``Robot as a service in cloud computing,'' in {\em 2010 Fifth IEEE International Symposium on Service Oriented System Engineering}, pp.~151--158, IEEE, 2010.

\bibitem{hu2012cloud}
G.~Hu, W.~P. Tay, and Y.~Wen, ``Cloud robotics: architecture, challenges and applications,'' {\em IEEE network}, vol.~26, no.~3, pp.~21--28, 2012.

\bibitem{kehoe2015survey}
B.~Kehoe, S.~Patil, P.~Abbeel, and K.~Goldberg, ``A survey of research on cloud robotics and automation,'' {\em IEEE Transactions on automation science and engineering}, vol.~12, no.~2, pp.~398--409, 2015.

\bibitem{buerkle2023towards}
A.~Buerkle, W.~Eaton, A.~Al-Yacoub, M.~Zimmer, P.~Kinnell, M.~Henshaw, M.~Coombes, W.-H. Chen, and N.~Lohse, ``Towards industrial robots as a service (iraas): Flexibility, usability, safety and business models,'' {\em Robotics and Computer-Integrated Manufacturing}, vol.~81, p.~102484, 2023.

\bibitem{straubinger2020overview}
A.~Straubinger, R.~Rothfeld, M.~Shamiyeh, K.-D. B{\"u}chter, J.~Kaiser, and K.~O. Pl{\"o}tner, ``An overview of current research and developments in urban air mobility--setting the scene for uam introduction,'' {\em Journal of Air Transport Management}, vol.~87, p.~101852, 2020.

\bibitem{cohen2021urban}
A.~P. Cohen, S.~A. Shaheen, and E.~M. Farrar, ``Urban air mobility: History, ecosystem, market potential, and challenges,'' {\em IEEE Transactions on Intelligent Transportation Systems}, vol.~22, no.~9, pp.~6074--6087, 2021.

\bibitem{dorling2016vehicle}
K.~Dorling, J.~Heinrichs, G.~G. Messier, and S.~Magierowski, ``Vehicle routing problems for drone delivery,'' {\em IEEE Transactions on Systems, Man, and Cybernetics: Systems}, vol.~47, no.~1, pp.~70--85, 2016.

\bibitem{scott2017drone}
J.~Scott and C.~Scott, ``Drone delivery models for healthcare,'' 2017.

\bibitem{dharmasena2019autonomous}
T.~Dharmasena, R.~de~Silva, N.~Abhayasingha, and P.~Abeygunawardhana, ``Autonomous cloud robotic system for smart agriculture,'' in {\em 2019 Moratuwa Engineering Research Conference (MERCon)}, pp.~388--393, IEEE, 2019.

\bibitem{milella2024robot}
A.~Milella, S.~Rilling, A.~Rana, R.~Galati, A.~Petitti, M.~Hoffmann, J.~L. Stanly, and G.~Reina, ``Robot-as-a-service as a new paradigm in precision farming,'' {\em IEEE access}, vol.~12, pp.~47942--47949, 2024.

\bibitem{yan2017cloud}
H.~Yan, Q.~Hua, Y.~Wang, W.~Wei, and M.~Imran, ``Cloud robotics in smart manufacturing environments: Challenges and countermeasures,'' {\em Computers \& Electrical Engineering}, vol.~63, pp.~56--65, 2017.

\bibitem{toscano2022integrating}
M.~Toscano-Moreno, J.~Bravo-Arrabal, M.~S{\'a}nchez-Montero, J.~S. Barba, R.~V{\'a}zquez-Mart{\'\i}n, J.~J. Fern{\'a}ndez-Lozano, A.~Mandow, and A.~Garcia-Cerezo, ``Integrating ros and android for rescuers in a cloud robotics architecture: Application to a casualty evacuation exercise,'' in {\em 2022 IEEE International Symposium on Safety, Security, and Rescue Robotics (SSRR)}, pp.~270--276, IEEE, 2022.

\bibitem{zhang2019data}
W.~Zhang, D.~Yang, and H.~Wang, ``Data-driven methods for predictive maintenance of industrial equipment: A survey,'' {\em IEEE systems journal}, vol.~13, no.~3, pp.~2213--2227, 2019.

\bibitem{carvalho2019systematic}
T.~P. Carvalho, F.~A. Soares, R.~Vita, R.~d.~P. Francisco, J.~P. Basto, and S.~G. Alcal{\'a}, ``A systematic literature review of machine learning methods applied to predictive maintenance,'' {\em Computers \& Industrial Engineering}, vol.~137, p.~106024, 2019.

\bibitem{wang2024dynamic}
L.~Wang, Z.~Zhu, and X.~Zhao, ``Dynamic predictive maintenance strategy for system remaining useful life prediction via deep learning ensemble method,'' {\em Reliability Engineering \& System Safety}, vol.~245, p.~110012, 2024.

\bibitem{ge_maintenance_insight}
{GE Aerospace}, ``Maintenance {Insight}.'' \url{https://www.geaerospace.com/systems/saas/maintenance-insight}, 2024.
\newblock Accessed on 2024-09-14.

\bibitem{mil-hdbk-217f}
{United States Department of Defense}, ``Military handbook: Reliability prediction of electronic equipment,'' Tech. Rep. MIL-HDBK-217F, Department of Defense, Washington, DC, USA, December 1991.

\bibitem{hinton2006prism}
A.~Hinton, M.~Kwiatkowska, G.~Norman, and D.~Parker, ``Prism: A tool for automatic verification of probabilistic systems,'' in {\em International conference on tools and algorithms for the construction and analysis of systems}, pp.~441--444, Springer, 2006.

\bibitem{vittinghoff2012regression}
E.~Vittinghoff, D.~V. Glidden, S.~C. Shiboski, and C.~E. McCulloch, {\em Regression Methods in Biostatistics: Linear, Logistic, Survival, and Repeated Measures Models}.
\newblock New York, NY: Springer, 2nd~ed., 2012.

\bibitem{karim2019reliability}
M.~R. Karim, M.~A. Islam, {\em et~al.}, {\em Reliability and survival analysis}.
\newblock Springer, 2019.

\bibitem{collett2015modelling}
D.~Collett, {\em Modelling Survival Data in Medical Research}.
\newblock Boca Raton, FL: Chapman and Hall/CRC, 3rd~ed., 2015.

\bibitem{eloranta2021cancer}
S.~Eloranta, K.~Smedby, P.~Dickman, and T.~Andersson, ``Cancer survival statistics for patients and healthcare professionals--a tutorial of real-world data analysis,'' {\em Journal of internal medicine}, vol.~289, no.~1, pp.~12--28, 2021.

\bibitem{papathanasiou2023machine}
D.~Papathanasiou, K.~Demertzis, and N.~Tziritas, ``Machine failure prediction using survival analysis,'' {\em Future Internet}, vol.~15, no.~5, p.~153, 2023.

\bibitem{cox1972regression}
D.~R. Cox, ``Regression models and life-tables,'' {\em Journal of the Royal Statistical Society: Series B (Methodological)}, vol.~34, no.~2, pp.~187--202, 1972.

\bibitem{kumar1994proportional}
D.~Kumar and B.~Klefsj{\"o}, ``Proportional hazards model: a review,'' {\em Reliability Engineering \& System Safety}, vol.~44, no.~2, pp.~177--188, 1994.

\bibitem{laffont2002theory}
J.-J. Laffont and D.~Martimort, {\em The theory of incentives: the principal-agent model}.
\newblock Princeton university press, 2002.

\bibitem{aswani2018inverse}
A.~Aswani, Z.-J. Shen, and A.~Siddiq, ``Inverse optimization with noisy data,'' {\em Operations Research}, vol.~66, no.~3, pp.~870--892, 2018.

\bibitem{chan2025inverse}
T.~C. Chan, R.~Mahmood, and I.~Y. Zhu, ``Inverse optimization: Theory and applications,'' {\em Operations Research}, vol.~73, no.~2, pp.~1046--1074, 2025.

\bibitem{aswani2019data}
A.~Aswani, Z.-J.~M. Shen, and A.~Siddiq, ``Data-driven incentive design in the medicare shared savings program,'' {\em Operations Research}, vol.~67, no.~4, pp.~1002--1026, 2019.

\bibitem{mintz2023behavioral}
Y.~Mintz, A.~Aswani, P.~Kaminsky, E.~Flowers, and Y.~Fukuoka, ``Behavioral analytics for myopic agents,'' {\em European journal of operational research}, vol.~310, no.~2, pp.~793--811, 2023.

\bibitem{li2025adaptive}
Q.~Li, K.~L. Gavin, C.~I. Voils, and Y.~Mintz, ``An adaptive optimization approach to personalized financial incentives in behavioral interventions,'' {\em Production and Operations Management}, p.~10591478251349391, 2025.

\bibitem{dogan2023estimating}
I.~Dogan, Z.-J.~M. Shen, and A.~Aswani, ``Estimating and incentivizing imperfect-knowledge agents with hidden rewards,'' {\em arXiv preprint arXiv:2308.06717}, 2023.

\bibitem{dogan2023repeated}
I.~Dogan, Z.-J.~M. Shen, and A.~Aswani, ``Repeated principal-agent games with unobserved agent rewards and perfect-knowledge agents,'' {\em arXiv preprint arXiv:2304.07407}, 2023.

\bibitem{scheid2024incentivized}
A.~Scheid, D.~Tiapkin, E.~Boursier, A.~Capitaine, E.~Moulines, M.~I. Jordan, E.-M. El-Mhamdi, and A.~Durmus, ``Incentivized learning in principal-agent bandit games,'' in {\em Proceedings of the 41st International Conference on Machine Learning}, ICML'24, JMLR.org, 2024.

\bibitem{agrawal2013thompson}
S.~Agrawal and N.~Goyal, ``Thompson sampling for contextual bandits with linear payoffs,'' in {\em International conference on machine learning}, pp.~127--135, PMLR, 2013.

\bibitem{dimakopoulou2019balanced}
M.~Dimakopoulou, Z.~Zhou, S.~Athey, and G.~Imbens, ``Balanced linear contextual bandits,'' in {\em Proceedings of the AAAI Conference on Artificial Intelligence}, pp.~3445--3453, 2019.

\bibitem{bastani2021mostly}
H.~Bastani, M.~Bayati, and K.~Khosravi, ``Mostly exploration-free algorithms for contextual bandits,'' {\em Management Science}, vol.~67, no.~3, pp.~1329--1349, 2021.

\bibitem{dann2022guarantees}
C.~Dann, Y.~Mansour, M.~Mohri, A.~Sekhari, and K.~Sridharan, ``Guarantees for epsilon-greedy reinforcement learning with function approximation,'' in {\em International conference on machine learning}, pp.~4666--4689, PMLR, 2022.

\bibitem{lobel2018multidimensional}
I.~Lobel, R.~P. Leme, and A.~Vladu, ``Multidimensional binary search for contextual decision-making,'' {\em Operations Research}, vol.~66, p.~1346–1361, Oct. 2018.

\bibitem{breslow1972contribution}
N.~E. Breslow, ``Contribution to discussion of paper by dr cox,'' {\em Journal of the Royal Statistical Society, Series B}, vol.~34, pp.~216--217, 1972.

\bibitem{Bertsekas2017VolI}
D.~P. Bertsekas, {\em Dynamic Programming and Optimal Control, Vol. I}.
\newblock Belmont, MA: Athena Scientific, 4th~ed., 2017.

\bibitem{Bertsekas2012VolII}
D.~P. Bertsekas, {\em Dynamic Programming and Optimal Control, Vol. II: Approximate Dynamic Programming}.
\newblock Belmont, MA: Athena Scientific, 4th~ed., 2012.

\end{thebibliography}
\bibliographystyle{ieeetr}

\end{document}